\documentclass[11pt,reqno,twoside]{amsart}
\usepackage{amsmath, amssymb, amscd, amsthm}
\newcommand{\R}{{\mathbb R}}
\newcommand{\N}{{\mathbb N}}
\newcommand{\C}{{\mathbb C}}
\newcommand{\Z}{{\mathbb Z}}
\newcommand{\Zp}{{\mathbb Z}_{\geq 0}}
\newcommand{\CR}{\C\backslash\R}
\newcommand{\E}{{\mathcal E}}
\newcommand{\hf}{{\frac12}}
\newcommand{\al}{\alpha}

\newcommand{\ga}{\gamma}

\newcommand{\de}{\delta}

\newcommand{\ep}{\varepsilon}

\newcommand{\te}{\theta}

\newcommand{\la}{\lambda}

\newcommand{\lt}{\ell^2(\Z)}
\newcommand{\vp}{\varphi}
\newcommand{\sgn}{{\rm sgn}}
\newcommand{\tvpo}[3]{ \,_2\vp_1 \left(
\genfrac{.}{.}{0pt}{}{#1}{#2} ;q,{#3}\right)}
\newcommand{\rvps}[5]{ \,_{#1}\vp_{#2}\left(
\genfrac{.}{.}{0pt}{}{#3}{#4} ;q,{#5}\right)}

\numberwithin{equation}{section}
\theoremstyle{plain}
\newtheorem{Thm}{Theorem}[section]
\newtheorem{Cor}[Thm]{Corollary}

\newtheorem{Lem}[Thm]{Lemma}
\newtheorem{Prp}[Thm]{Proposition}

\theoremstyle{remark}
\newtheorem{Rmk}[Thm]{Remark}

\begin{document}
\title[Orthogonality relations]{One-parameter orthogonality relations \\ for
basic hypergeometric series}
\author{Erik Koelink}
\address{Technische Universiteit Delft,
Faculteit Elektrotechniek, Wiskunde en Informatica, TWA, Postbus 5031,
2600 GA Delft, the Netherlands.}
\email{h.t.koelink@math.tudelft.nl}

\dedicatory{Dedicated to Tom Koornwinder on
the occasion of his 60th birthday}

\date{May 27, 2003
\\ \indent 2000 {\it Mathematics Subject Classification}.
33D45, 33D15.}

\keywords{basic hypergeometric series, spectral analysis,
orthogonality relations}

\begin{abstract} The second order hypergeometric $q$-difference
operator is studied for the value $c=-q$. For certain parameter
regimes the corresponding recurrence relation can be related to
a symmetric operator on the Hilbert space $\lt$. The operator
has deficiency indices $(1,1)$ and we describe as explicitly
as possible the spectral resolutions of the self-adjoint
extensions. This gives rise to one-parameter orthogonality
relations for sums of two ${}_2\vp_1$-series. In particular, we find
that the Ismail-Zhang $q$-analogue of the exponential function
satisfies certain orthogonality relations.
\end{abstract}

\maketitle


\section{Introduction}\label{sec:intro}

As is well known, special functions arise in several contexts
in mathematics. One of the areas is the theory of self-adjoint
operators on a Hilbert space, see e.g. Titchmarsh \cite{Titc}.
On the one hand, given an explicit self-adjoint operator,
we can try to use special functions in order to obtain the
spectral decomposition of the self-adjoint operator. On the
other hand, given a family of interesting special functions, 
we can try
to find a self-adjoint operator which has these interesting
special functions as eigenfunctions. If we can give the
spectral decomposition of the corresponding operator we can use
this information to find e.g. orthogonality relations
or a corresponding integral transform for
the special function we have started with.

In this paper we are in the second situation. The interesting
special function is
\begin{equation}\label{eq:defcurlyE}
\E_q(z;t) = \frac{(-t;q^\hf)_\infty}{(qt^2;q^2)_\infty}
\,_2\vp_1 \left(
\genfrac{.}{.}{0pt}{}{q^{\frac14}y,q^{\frac14}/y}{-q^\hf} ;q^\hf,-t\right)
\end{equation}
originally introduced by Ismail and Zhang
\cite[(1.22) with $a=i$, $b=-it$]{IsmaZ} up to
a normalisation factor,
and the expression used here can be found in
Ismail and Stanton \cite[Corr.~4.3]{IsmaSCJM},
\cite[Corr.~2.5]{IsmaSpp}, see also Suslov \cite{Susl} for more
information. Ismail and Zhang
\cite[(1.25)]{IsmaZ}, see also \cite{IsmaSCJM}, point out that formally
\begin{equation}\label{eq:limcurlyE}
\lim_{q\uparrow 1}\
\E_q(z;\hf t)\Big\vert_{t=(1-q)\la}
= e^{\la z}.
\end{equation}
This function has been studied intensively recently because
it is the appropriate $q$-analogue of the exponential function
well suited for the Askey-Wilson difference operator.

In the point of view of this paper, we study the function
$\E_q(z,t)$ as a function of $t$. The parameter $z$ occurs
as the spectral parameter. Because of the expression
\eqref{eq:defcurlyE} $\E_{q^2}$ is an eigenfunction of the
second order hypergeometric $q$-difference operator. It
is convenient to switch to $q^2$ in \eqref{eq:defcurlyE}.
This operator and its eigenfunctions have
been studied in connection with
representation theory of non-compact quantum groups, in
particular the quantum analogue of $SU(1,1)$,
see \cite{Kake}, \cite{KoelSRIMS}, and \cite{KoelSNATO}
for a more general scheme. The parameter regimes for the
basic hypergeometric function in these papers does not
include the case corresponding to $\E_q$, so we have to perform the
spectral analysis again. The crucial property is that the
lower parameter $c=-q$.

It turns out that for specific values of the remaining
parameters the second order hypergeometric $q$-difference operator
can be realised as an unbounded symmetric operator on the Hilbert
space $\lt$ of square integrable sequences. In particular,
this occurs for $\E_{q^2}$. However, it turns out that the
corresponding operator is not essentially self-adjoint,
but it has deficiency indices $(1,1)$. We describe
the self-adjoint extensions, which depend on one extra
parameter, and we study the corresponding spectral
decompositions. There is always continuous spectrum
on $[-1,1]$, and the point spectrum is an infinite set
tending to plus and/or minus infinity which is described
as the zero set of some explicit function. For the case
of the function $\E_{q^2}$ we establish that this set
consists of two $q^2$-quadratic grids. The corresponding
transforms do not give orthogonality relations for
$\E_{q^2}$ but for a linear combination of two $\E_{q^2}$'s
similar to the relation $2\cos \la x = e^{i\la x}+e^{-i\la x}$.
So we can think of the result as a $q$-analogue of the
Fourier-cosine transform instead of the Fourier transform.
We perform the spectral analysis in somewhat greater
generality, and the main result is Theorem
\ref{thm:case1spectraldecomp} and its counterpart
Theorem \ref{thm:case2spectraldecomp} for
another parameter regime.

We also present the link with the recurrence relation
for the big $q$-Jacobi functions \cite{KoelSCA} or
the associated dual $q$-Hahn polynomials \cite{GuptIM},
and this leads to a quadratic transformation in which
a ${}_2\vp_1$-series in base $q$ is given as a
${}_3\vp_2$-series in base $q^2$. In particular, this
gives a new expression for $\E_q$ as a sum of
two ${}_3\vp_2$-series in base $q$.

The plan of the paper is as follows. In
\S \ref{sec:diffeqsols} we recall the second order
hypergeometric $q$-difference equation, its solutions and
their interrelations. In \S \ref{sec:symmformdiffeq} we
discuss for which parameter regimes the recurrence relation
can be interpretated as a symmetric operator on $\lt$.
We recall the general theory of doubly infinite Jacobi
operators on $\lt$ in \S \ref{sec:gendJacobiop}. In
\S\S \ref{sec:specdecL1}, \ref{sec:specdecL2}
we work out the spectral decompositions of the
self-adjoint extensions as explicitly as possible.
In \S  \ref{sec:specdecL1} we give detailed arguments,
and we indicate the (similar) arguments in the easier
case of \S \ref{sec:specdecL2}. Finally, in
\S \ref{sec:quadrtrans} we indicate the link with
the big $q$-Jacobi functions, and we derive the
quadratic relation.

{\em Notation.}
In this paper we follow the notation for basic hypergeometric
series of Gasper and Rahman \cite{GaspR}. Our standing
assumption on $q$ is $0<q<1$. The series
\begin{equation}\label{eq:defbhs}
\rvps{r+1}{r}{a_1,\ldots,a_{r+1}}{b_1,\ldots,b_r}{z}
=\sum_{k=0}^\infty \frac{(a_1,\ldots,a_{r+1};q)_k}
{(b_1,\ldots,b_r,q;q)_k} z^k,
\end{equation}
where $(a;q)_k=\prod_{i=0}^{k-1}(1-aq^i)$,
$k\in\Zp\cup\{\infty\}$, and
$(a_1,\ldots,a_r;q)_k=\prod_{i=1}^k (a_i;q)_k$.
Generically the radius of convergence of the series in \eqref{eq:defbhs}
is $1$, but the series has a unique analytic continuation
to $\C\backslash[1,\infty)$.
We also use $\te(z)=(z,q/z;q)_\infty$ for the (renormalised)
Jacobi theta-function, and
$\te(a_1,\ldots,a_r)=\te(a_1)\cdots \te(a_r)$. The identity
\begin{equation}\label{eq:thetaidentity}
\te(aq^k) = (-a)^{-k} q^{-\hf k(k-1)}\, \te(a)
\end{equation}
is useful.

{\em Dedication and acknowledgement.}
The paper is dedicated to Tom Koornwinder from whom I have
learnt very much about special functions and
representation theory. In particular, his unpublished
notes on spectral theory for ${}_2\vp_1$-series have been
very influential for this paper.
I thank Mourad Ismail for discussions about the
$q$-analogue of the exponential function defined in
\eqref{eq:defcurlyE} that has triggered the research
for this paper, for his kind hospitality during
a visit when these discussions took place, and for
his comments on a previous version of this paper. 


\section{The difference equation and its solutions}
\label{sec:diffeqsols}

In this section we consider the second order
hypergeometric $q$-difference operator for which
the ${}_2\vp_1$-series in the definition
\eqref{eq:defcurlyE} is an eigenfunction. We
discuss other solutions and their connection
coefficients.

The hypergeometric difference equation is
\begin{equation}\label{eq:hypergeomdiffop}
(c-abz)\, f(qz) + \bigl(-(c+q)+(a+b)z\bigr)\, f(z)
+(q-z)\, f(z/q)=0,
\end{equation}
see \cite[Exerc. 1.13]{GaspR}, having
$$
f(z) = \tvpo{a,b}{c}{z}
$$
as a solution.
We are particularly
interested in the case $c=-q$, cf. \eqref{eq:defcurlyE}.

\begin{Lem}\label{lem:diffeqsol} The difference equation
\begin{equation}\label{eq:diffequnnorm}
2z\, f_k(z) = \frac{1+a^2tq^{k-1}}{atq^{k-1}}\, f_{k+1}(z)
- \frac{1-q^{k-1}t}{atq^{k-1}}\, f_{k-1}(z), \qquad k\in\Z,
\end{equation}
is solved by, where $z=\hf(y+y^{-1})$,
\begin{equation*}
\begin{split}
u_k(z) &= \tvpo{ay,a/y}{-q}{tq^k}, \\
v_k(z) &= (-1)^k \tvpo{-ay,-a/y}{-q}{tq^k}, \\
F_k(y) &= (ay)^{-k} \tvpo{ay,-ay}{qy^2}{-\frac{q^{2-k}}{a^2t}},
\qquad y^2 \notin q^{-\N}.
\end{split}
\end{equation*}
\end{Lem}

Here, and in the sequel, we always assume that $a\not=0$,
$t\not=0$.

\begin{proof} This is a straightforward verification
using \eqref{eq:hypergeomdiffop}.
\end{proof}

In \S \ref{sec:quadrtrans} we also give expressions for the solutions of
\eqref{eq:diffequnnorm} in terms of ${}_3\vp_2$-series
using a quadratic transformation.

\begin{Rmk}\label{rmk:symmdifeq}
Note that the difference equation \eqref{eq:diffequnnorm}
has two obvious symmetries. The
first is $a\leftrightarrow -a$, $z\leftrightarrow -z$,
leaving all solutions unchanged. The second
symmetry is $a\leftrightarrow -a$, $f_k\leftrightarrow (-1)^kf_k$,
which interchanges the solutions $u_k\leftrightarrow v_k$ and
leaves $F_k$ unchanged.
\end{Rmk}

\begin{Rmk} We are mainly interested in the case
that the coefficients in \eqref{eq:diffequnnorm} do not vanish
for $k\in\Z$, i.e. we assume $t\notin q^\Z$, $-ta^2\notin q^\Z$.
In case one of the coefficients does vanish, we can assume
without loss of generality that $t=q$ or $a^2t=-q$. In this
case the recurrence can be split into two recurrence relations
labeled by $\N$. So we have polynomial solutions, and the
polynomials can be given explicitly in terms of symmetric
Al-Salam--Chihara polynomials in base $q$ for negative $k$
and in base $q^{-1}$ for positive $k$.
\end{Rmk}

Since Lemma \ref{lem:diffeqsol} describes four solutions
(note that $F_k(y^{-1})$ is also a solution) to
\eqref{eq:hypergeomdiffop}, whose solution space is
two-dimensional, we find relations between the
solutions.

\begin{Lem}\label{lem:relsol}
The solutions of Lemma \ref{lem:diffeqsol} are related by,
$z=\hf(y+y^{-1})$,
\begin{equation}\label{eq:cexp}
\begin{split}
u_k(z) &= c(y;a,t) \, F_k(y) + c(y^{-1};a,t)\, F_k(y^{-1}), \\
v_k(z) &= c(y;-a,t) \, F_k(y) + c(y^{-1};-a,t)\, F_k(y^{-1}), \\
c(y;a,t) &=
\frac{(a/y,-q/ay,ayt,q/ayt;q)_\infty}{(-q,y^{-2},t,q/t;q)_\infty},
\end{split}
\end{equation}
and
\begin{equation}\label{eq:dexp}
\begin{split}
F_k(y) &= d(y;a,t) \, u_k(z) + d(y;-a,t)\, v_k(z), \\
d(y;a,t) &=
\frac{(-ay,qy/a,-at/qy,-q^2y/at;q)_\infty}
{(-1,qy^2,-a^2t/q,-q^2/a^2t;q)_\infty},
\end{split}
\end{equation}
\end{Lem}

\begin{proof} The second equation
of \eqref{eq:cexp} follows from the first
using the symmetries as in Remark \ref{rmk:symmdifeq}.
The first equation of \eqref{eq:cexp} follows from
\cite[(4.3.2)]{GaspR}. The expansion \eqref{eq:dexp}
can be proved similarly, or by inverting \eqref{eq:cexp}.
In the last case the addition formula for Jacobi theta
functions has to be used, see e.g. \cite[Exerc.~2.16]{GaspR}, to
find
\begin{equation*}
\det\left(\begin{matrix} c(y;a,t) & c(y^{-1};a,t) \\
c(y;-a,t) & c(y^{-1};-a,t)
\end{matrix}\right) = \frac{-aty \, \te(-a^2t,t^{-1},y^{-2},-1)}
{(-q,-q,y^2,y^{-2},t,q/t,t,q/t;q)_\infty}
= \frac{2a}{y^{-1}-y}\frac{\te(-a^2t)}{\te(t)}.
\end{equation*}
\end{proof}


\section{Symmetric form of the difference equation}
\label{sec:symmformdiffeq}

Since we want to find a symmetric operator on
the Hilbert space $\lt$ for
which the $\E_{q^2}(z;tq^k)$ occur as eigenfunctions
we need to find conditions on $a$ and $t$ such that
we can rewrite \eqref{eq:diffequnnorm} in
a symmetric form. This is done in this section.

Let $f_k(z)$ satisfy \eqref{eq:diffequnnorm}, then
$g_k(z)=\al_k f_k(z)$, for non-zero constants $\al_k$, satisfies
\begin{equation}\label{eq:symmform1}
2z\, g_k(z) = \frac{\al_k}{\al_{k+1}}
\frac{1+a^2tq^{k-1}}{atq^{k-1}}\, g_{k+1}(z)
- \frac{\al_k}{\al_{k-1}}
\frac{1-q^{k-1}t}{atq^{k-1}}\, g_{k-1}(z).
\end{equation}
We need to determine if
we can rewrite the recurrence in the symmetric form
\begin{equation}\label{eq:symmform}
2z\, g_k(z) = a_k\, g_{k+1}(z) +  a_{k-1}\, g_{k-1}(z),
\qquad a_k>0.
\end{equation}
From the coefficient of $g_{k+1}(z)$ in equations
\eqref{eq:symmform1} and \eqref{eq:symmform}
we find the first equality
$$
\bar a_{k-1} =  \frac{\bar\al_{k-2}}{\bar\al_k}
\frac{1+\bar a^2\bar tq^{k-2}}{\bar a\bar tq^{k-1}} =
- \frac{\al_k}{\al_{k-1}}
\frac{1-q^{k-1}t}{atq^{k-1}},
$$
where the second equality follows from the coefficient for
$g_{k-1}(z)$ in \eqref{eq:symmform1} and \eqref{eq:symmform}.
Hence, \begin{equation}\label{eq:symmcondat}
\left| \frac{\al_{k-1}}{\al_k}\right|^2 =
- \frac{1-q^{k-1}t}{1+\bar a^2\bar t q^{k-2}}
\, \frac{\bar a\bar t}{qat}
\end{equation}
and we can make the appropriate choice for $\al_k$
precisely when
the right hand side of \eqref{eq:symmcondat} is strictly
positive. Note that the choice for $\al_k$ is determined
by \eqref{eq:symmcondat} and one initial value, say for $\al_0$,
up to a phase factor. We can choose the phase factor such that
the value $a_k$ in \eqref{eq:symmform} is indeed positive.

\begin{Lem}\label{lem:symmform} In the following cases \eqref{eq:symmform1}
can be written in the symmetric form \eqref{eq:symmform}
with $a_k>0$, $\forall k\in\Z$:
\begin{enumerate}
\item $a=\sqrt{q}\, e^{i\psi}$, $t=ire^{-i\psi}$ with
$r\in\R\backslash\{0\}$,
\item $t<0$, $a=is$ with $s\in\R\backslash\{0\}$,
\item $t\in \R_{>0}$, $a=is$, $s\in\R\backslash\{0\}$, such that
there exists $k_0\in\Z$ with $tq^{k_0+1}<1<tq^{k_0}$ and
$s^2tq^{k_0}<1<s^2tq^{k_0-1}$.
\end{enumerate}
\end{Lem}

Note that there is overlap between cases (1) and (2), and  (1)
and (3).
For the remainder of the paper we stick to the cases (1) and (2),
where in case (1) we moreover assume that $t\notin\R_{>0}$ in
order to have the ${}_2\vp_1$-series in Lemma \ref{lem:diffeqsol}
well-defined as analytic functions on $\C\backslash[1,\infty)$.

We fix the corresponding values of the coefficients $a_k$
and $\al_k$ as follows. In case (1) we take
\begin{equation}\label{eq:ak1}
\begin{split}
a_k &=|r|^{-1}q^{-k}\sqrt{1-2rq^k\sin\psi+r^2q^{2k}} =
\left|\frac{1+ire^{i\psi}q^k}{irq^k}\right|, \\
\al_k &= e^{i\phi_k}q^{\hf k}, \qquad
\phi_{k+1}-\phi_k \equiv \arg (1+ire^{i\psi}q^k)-\hf\pi\, \sgn(r) \mod 2\pi ,
\end{split}
\end{equation}
and in case (2) we take
\begin{equation}\label{eq:ak2}
\begin{split}
a_k & = \frac{q^{\hf-k}}{|st|}
\sqrt{(1-tq^k)(1-ts^2q^{k-1})}
= \sqrt{(1-q^{-k}/t)(1-q^{1-k}/ts^2)} , \\
\al_k &= i^k\, s^k
\sqrt{\frac{(q^{2-k}/s^2t;q)_\infty}{(q^{1-k}/t;q)_\infty}}
=
(i\, \sgn(s))^k q^{\hf k}
\sqrt{\frac{(tq^k;q)_\infty\,  \te(s^2t/q)}
{(ts^2q^{k-1};q)_\infty\, \te(t)}},
\end{split}
\end{equation}
using the $\te$-product identity \eqref{eq:thetaidentity}.

It follows from Lemma \ref{lem:relsol} that in all cases
of Lemma \ref{lem:symmform} we have
\begin{equation}\label{eq:relforcf}
\overline{c(\bar y;a,t)} =  \frac{\te(t)}{\te(\bar t)}\,
c(y;-a,t).
\end{equation}


\section{Generalities on doubly infinite Jacobi operators}
\label{sec:gendJacobiop}

In this section we recall some of the general theory for
the spectral analysis of doubly infinite Jacobi operators
on the Hilbert space $\lt$ given by \eqref{eq:defLonlt}.
In the cases considered in this paper we have to deal with
one-dimensional deficiency spaces, and the self-adjoint extensions are
described. The results of this section can be found in 
\cite{MassR}, \cite{KoelLaredo}, \cite[Ch. 7]{Bere}, 
\cite{Prui} and for more generalities
Dunford and Schwartz \cite{DunfS} can be consulted.

We consider next the corresponding operator $L$ on
the Hilbert space $\lt$ equipped with an orthonormal
basis $\{ e_k\}_{k\in\Z}$ defined by
\begin{equation}\label{eq:defLonlt}
2 L\, e_k = a_k\, e_{k+1} + a_{k-1}\, e_{k-1}
\end{equation}
with $a_k>0$ as in \eqref{eq:ak1} and \eqref{eq:ak2}.

The operator is initially defined on the dense
domain ${\mathcal D}$ of finite linear combinations of
the basis vectors $e_k$. The operator $(L, {\mathcal D})$ is
a symmetric operator, and,
since $a_k\in\R$, $L$ commutes with conjugation. So the
deficiency indices are equal, and since the solution space
of $L\, \xi=z\, \xi$ is two-dimensional the deficiency indices
are $(0,0)$, $(1,1)$, or $(2,2)$. In the cases
\eqref{eq:ak1} and \eqref{eq:ak2}
it follows that $a_k$ is bounded for $k\to-\infty$.
By Theorem~2.1 of Masson and Repka \cite{MassR},
see also \cite[(4.2.2)]{KoelLaredo},
we find that deficiency indices are $(0,0)$ or $(1,1)$.
The adjoint operator is $(L^\ast, {\mathcal D}^\ast)$ given by
\begin{equation*}
L^\ast \Bigl(\sum_{k=-\infty}^\infty \xi_k\, e_k\Bigr)
= \sum_{k=-\infty}^\infty (a_k \xi_{k+1} +
a_{k-1}\xi_{k-1}) \,  e_k, \qquad
{\mathcal D}^\ast = \bigl\{ \xi\in \lt \mid
L^\ast \xi\in\lt\bigr\}.
\end{equation*}
In the cases considered in this paper the deficiency indices
are $(1,1)$.

Note that $g(z) = \sum_{k=-\infty}^\infty
g_k(z)\, e_k$ is a solution to the eigenvalue equation
$L^\ast \xi=z\, \xi$ precisely when $g_k(z)$ satisfies
the recurrence relation \eqref{eq:symmform}. We
denote by $\al f(z)$ the solutions to the eigenvalue
equation of the form
$\al f(z) = \sum_{k=-\infty}^\infty \al_k f_k(z)\, e_k$
with $f_k(z)$ a solution to the recurrence of
Lemma \ref{lem:diffeqsol} and $\al_k$ as in
\eqref{eq:ak1} or \eqref{eq:ak2}.

Recall the Wronskian (or Casorati determinant),
\begin{equation}\label{eq:defWronskian}
[u,v]_k = \hf a_k\, (u_{k+1}\, v_k-u_k\, v_{k+1}),
\qquad u=\sum_{k=-\infty}^\infty u_k\, e_k, \quad
v=\sum_{k=-\infty}^\infty v_k\, e_k.
\end{equation}
If moreover $u$ and $v$ satisfy the eigenvalue equation
$L^\ast\, \xi=z\, \xi$, then $[u,v]_k$ is independent
of $k\in\Z$. And $u$ and $v$ are linearly independent
solutions of the eigenvalue equation if and only
if the Wronskian $[u,v]\not=0$. Note that we do not
impose $u,v\in\lt$.

Since $a_k$ is bounded as $k\to-\infty$ the space
$$
S^-(z) = \Bigl\{ \xi=\sum_{n=-\infty}^\infty \xi_n\, e_n
\mid L^\ast\, \xi=z\, \xi,\ \sum_{n=-\infty}^0 |\xi_n|^2<\infty
\Bigr\}
$$
is one-dimensional for $z\in\CR$.
We assume it is spanned by
$\Psi(z)=\sum_{k=-\infty}^\infty \Psi_k(z)\, e_k$
satisfying $\Psi_k(z)= \overline{\Psi_k(\bar z)}$.
Note that this condition can be imposed since $L$ commutes
with complex conjugation.
The similarly defined space
$$
S^+(z) = \Bigl\{ \xi=\sum_{n=-\infty}^\infty \xi_n\, e_n
\mid L^\ast\, \xi=z\, \xi,\ \sum_{n=0}^\infty |\xi_n|^2<\infty\Bigr\}
$$
is at most two-dimensional and at least one-dimensional
for $z\in\CR$.
We show later that in cases (1) and
(2) of Lemma \ref{lem:symmform}
the space $S^+(z)$ is two-dimensional, so that
the deficiency indices of $(L,{\mathcal D})$
are $(1,1)$. Indeed, $\dim \ker(L^\ast \pm i)=
\dim S^+(\mp i)\cap S^-(\mp i)= 1$.
The fact $\dim S^+(z)=2$, $z\in\CR$,
follows from the fact that the
asymptotic behaviour of $u_k(z)$ is the same as that
of $v_k(z)$ (up to a sign $(-1)^k$) as $k\to\infty$, assuming
we know that the solutions $u_k(z)$ and $v_k(z)$
are linearly independent,
see \S\S \ref{sec:specdecL1}, \ref{sec:specdecL2}.
So we have
$\dim \ker(L^\ast \pm i)=1$, and $\Psi(\pm i)\in
\ker(L^\ast\mp i)$. Then the self-adjoint extensions of
$(L,{\mathcal D})$ are given by $(L^\ast, {\mathcal D}_\te)$
with
\begin{equation}\label{eq:defDtheta}
{\mathcal D}_\te = \Bigl\{ \xi\in {\mathcal D}^\ast \mid
\lim_{N\to\infty} [\xi, e^{i\te}\Psi(i) + e^{-i\te}\Psi(-i)]_N = 0
\Bigr\}, \quad \te\in[0,2\pi).
\end{equation}

Pick $\overline{\psi(\bar z)}\in S^+(z)\cap {\mathcal D}_\te$,
then we can
describe the resolvent for the self-adjoint operator
$(L,{\mathcal D}_\te)$ in terms of the Green function
\begin{equation}\label{eq:defGreen}
G_{k,l}(z) = \frac{1}{[\Psi(z), \overline{\psi(\bar z)}]}
\begin{cases}
\Psi_k(z)\overline{\psi_l(\bar z)}, & k\leq l, \\[2pt]
\Psi_l(z)\overline{\psi_k(\bar z)}, & l\leq k,
\end{cases}
\end{equation}
and the resolvent $R(z)=(z-L)^{-1}$ is given by, $e\in\lt$,
$$
R(z)\xi = \sum_{k=-\infty}^\infty \bigl(R(z)\xi\bigr)_k\, e_k, \qquad
\bigl(R(z)\xi\bigr)_k = \sum_{l=-\infty}^\infty \xi_l G_{k,l}(z)
=\langle \xi, \overline{G_{k,\cdot}(z)}\rangle.
$$
Note that for $\xi,\eta\in\lt$
\begin{equation}\label{eq:exprresolvinGreen}
[\Psi(z), \overline{\psi(\bar z)}]\, \langle R(z)\xi,\eta\rangle
=\sum_{k\leq l} \Psi_k(z) \overline{\psi_l(\bar z)}
\bigl( \xi_l\bar \eta_k + \xi_k \bar \eta_l \bigr) (1-\hf \de_{k,l}).
\end{equation}
The corresponding spectral measure $E$ of the self-adjoint
operator $(L,{\mathcal D}_\te)$ can be obtained from
the resolvent by
\begin{equation}\label{eq:StieltjesPerron}
E_{\xi,\eta}\bigl( (x_1,x_2)\bigr) = \lim_{\de\downarrow 0}
\lim_{\ep\downarrow 0} \frac{1}{2\pi i} \int_{x_1+\de}^{x_2-\de}
\langle R(x-i\ep)\xi,\eta\rangle - \langle R(x+i\ep)\xi,\eta\rangle \, dx
\end{equation}
for $\xi,\eta\in\lt$.


\section{Spectral decomposition of $L$ in the first case}
\label{sec:specdecL1}

In this section we calculate the spectral measure as
explicitly as possible of the self-adjoint extensions of
$(L,{\mathcal D})$ with $L$ as in \eqref{eq:defLonlt} with $a_k$
given by \eqref{eq:ak1}. This depends on the parameter $\te$
of the self-adjoint extension $(L^\ast,{\mathcal D}_\te)$ of
$(L,{\mathcal D})$. There is always continuous spectrum on the
interval $[-1,1]$, and an infinite series of discrete mass
points tending to plus or minus $\infty$. The location of the
discrete mass points depends on the choice of the
self-adjoint extension. In this section we always
have $a=q^\hf e^{i\psi}$ and $t=ire^{-i\psi}$ as
in case (1) of Lemma \ref{lem:symmform}, but we
keep the notation $a$ and $t$ in order to keep the analogy with
\S \ref{sec:specdecL1} in \S \ref{sec:specdecL2}.

Using \eqref{eq:ak1} we see that
\begin{equation}\label{eq:WFFcaseone1}
a_k e^{i(\phi_{k+1}-\phi_k)} = \frac{1+ire^{i\psi}q^k}{irq^k}.
\end{equation}

\begin{Lem}\label{lem:conjFkyreltoFky}
There is a $\ga\in\R$ such that
$\overline{e^{i\ga}\al_k F_k(\bar y)} =
e^{i\ga}\al_k F_k( y)$.
\end{Lem}

\begin{proof}
Now
\begin{equation}\label{eq:conjFkyreltoFky}
\begin{split}
\overline{\al_k F_k(\bar y)} &=
\overline{\al_k} (q^\hf e^{-i\psi}y)^{-k}
\tvpo{q^\hf e^{-i\psi}y,-q^\hf e^{-i\psi}y}{qy^2}{-ie^{-i\psi}
\frac{q^{1-k}}{r}} \\
&= e^{2ik\psi} \frac{(ie^{-i\psi}q^{1-k}/r;q)_\infty}
{(-ie^{i\psi}q^{1-k}/r;q)_\infty} \overline{\al_k} F_k(y) \\
&= e^{2ik\psi-2i\phi_k} \frac{(ie^{-i\psi}q^{1-k}/r;q)_\infty}
{(-ie^{i\psi}q^{1-k}/r;q)_\infty} \, \al_k F_k(y)
\end{split}
\end{equation}
where we use \cite[(1.4.6)]{GaspR} in the second equation.
From this calculation we only find
$\overline{\al_k F_k(\bar y)} = C_k \al_k F_k(y)$
with $|C_k|=1$. It remains to show that $C_k$ is independent of $k$,
and this follows from
\begin{equation*}
\begin{split}
\frac{C_{k+1}}{C_k} & = e^{2i(\phi_k-\phi_{k+1})} e^{2i\psi}
\frac{1-ie^{-i\psi}q^{-k}/r}{1+ie^{i\psi}q^{-k}/r} \\
&= \left( e^{i(\phi_k-\phi_{k+1})} \frac{1+irq^ke^{i\psi}}{irq^k}
\right) \left( e^{i(\phi_k-\phi_{k+1})}
\frac{-irq^k}{1-ire^{-i\psi}q^k}\right) =
\frac{a_k}{\overline{a_k}} =1
\end{split}
\end{equation*}
by \eqref{eq:WFFcaseone1} and $a_k\in\R$.
\end{proof}

\begin{Rmk} Note that the one-dimensional space
$S^-(z)$, $z\in\CR$, is spanned by $\al F(y)$, $|y|<1$,
with $z=\hf(y+y^{-1})$, since
$|\al_k F_k(y)| = {\mathcal O}(|y|^{-k})$
as $k\to-\infty$.
Since the coefficients $a_k$ are positive
$S^-(z)$ is also spanned $\overline{\al F(\bar y)}$, so we
see that $\overline{\al F(\bar y)} = C \al F(y)$ for
some constant $C$. In Lemma \ref{lem:conjFkyreltoFky}
we have shown moreover that $C$ is independent of $z$.
\end{Rmk}

A straightforward corollary to
Lemma \ref{lem:conjFkyreltoFky} is
\begin{equation}\label{eq:conjukzreltovkz}
\overline{e^{i\ga}\te(t) \al_k u_k(\bar z)} =
\te(t) e^{i\ga} \al_k v_k(z)
\end{equation}
using Lemma \ref{lem:relsol} and \eqref{eq:relforcf}.

It follows from
Lemma \ref{lem:conjFkyreltoFky}
that, in the notation of \S \ref{sec:gendJacobiop},
we have  $\Psi(z)=e^{i\ga}\al F( y) =
\overline{e^{i\ga} \al F(\bar y)}$
with $z=\hf(y+y^{-1})$ and $|y|<1$.

As is clear from \S \ref{sec:gendJacobiop} we need to
calculate various Wronskians in order to determine
the domain of the self-adjoint extensions and
the corresponding spectral measures. We state the
results in the following lemma.

\begin{Lem}\label{lem:allWronskians} We have the following Wronskians;
\begin{equation*}
\begin{split}
[\al F(y), \overline{\al F(\bar y^{-1})}] &= \hf (y^{-1}-y), \\
\lim_{N\to\infty}
[\overline{\al u(\bar w)}, \al F(y)]_N  &=
-\frac{q^\hf}{ir} d(y;a,t), \\
\lim_{N\to\infty}
[\overline{\al v(\bar w)}, \al F(y)]_N &=
\frac{q^\hf}{ir} d(y;-a,t)
\end{split}
\end{equation*}
\end{Lem}

\begin{proof} We first calculate the Wronskian
\begin{equation*}
\begin{split}
[\al F(y),  \overline{\al F(\bar y^{-1})}]_k =&
\frac{a_k}2 \bigl( \al_{k+1} F_{k+1}(y) \, \overline{\al_k
 F_k(\bar y^{-1})} -
\al_k F_k(y)\,  \overline{\al_{k+1} F_{k+1}(\bar y^{-1})}\bigr) \\
=& \hf q^{k+\hf} \bigl( a_k e^{i(\phi_{k+1}-\phi_k)}
\, F_{k+1}(y) \overline{F_k(\bar y^{-1})} -
a_k e^{i(\phi_k-\phi_{k+1})}
\, F_k(y) \overline{F_{k+1}(\bar y^{-1})}\bigr) \\
=& \frac{q^{k+\hf}}2
\frac{(1+ire^{i\psi}q^k)}{irq^k}  \bigl( q^\hf
e^{i\psi}y\bigr)^{-k-1} \overline{\bigl( q^\hf
e^{i\psi}\bar y^{-1}\bigr)^{-k}}(1+{\mathcal O}(q^{-k})) \\
&\qquad -
\frac{q^{k+\hf}}2
\frac{(1-ire^{-i\psi}q^k)}{-irq^k}
\bigl( q^\hf e^{i\psi}y\bigr)^{-k} \overline{\bigl( q^\hf
e^{i\psi}\bar y^{-1}\bigr)^{-k-1}}(1+{\mathcal O}(q^{-k})) \\
=& \hf (y^{-1}-y) \bigl(1+{\mathcal O}(q^{-k})\bigr)
\end{split}
\end{equation*}
using \eqref{eq:WFFcaseone1} twice and the
expression for $F_k(y)$ in Lemma \ref{lem:diffeqsol}. Since
the Wronskian is independent of $k$, we let $k\to-\infty$
to find the first statement of the lemma.

The next statement of the lemma follows from
\begin{equation*}
\begin{split}
&\lim_{N\to\infty}
[\overline{\al u(\bar w)}, \al F(y)]_N  \\= &
\lim_{N\to \infty}
\hf q^{N+\hf} \Bigl( a_Ne^{i\phi_N-i\phi_{N+1}}
\overline{u_{N+1}(\bar w)}
\bigl\{ d(y;a,t)u_N(z) + d(y;-a,t) v_N(z)\bigr\} \\
&\qquad\qquad\qquad - a_Ne^{i\phi_{N+1}-i\phi_N}
\overline{u_N(\bar w)}
\bigl\{d(y;a,t) u_{N+1}(z) + d(y;-a,t)v_{N+1}(z)\bigr\}\Bigr) \\
= & \lim_{N\to \infty} \hf q^{N+\hf} \Bigl(
\frac{1-ire^{-i\psi}q^N}{-irq^N}(1+{\mathcal O}(q^N))
\bigl\{ d(y;a,t)(1+{\mathcal O}(q^N)) +
d(y;-a,t)(-1)^N(1+{\mathcal O}(q^N))\bigr\}\\
&\qquad\quad\quad
 - \frac{1+ire^{i\psi}q^N}{irq^N}(1+{\mathcal O}(q^N))
\bigl\{ d(y;a,t)(1+{\mathcal O}(q^N)) +
 d(y;-a,t)(-1)^{N+1}(1+{\mathcal O}(q^N))\bigr\}\Bigr) \\
= & -\frac{q^\hf}{ir} d(y;a,t)
\end{split}
\end{equation*}
using \eqref{eq:WFFcaseone1} and Lemma
\ref{lem:diffeqsol}. The last statement follows
similarly.
\end{proof}

It follows from Lemma \ref{lem:allWronskians}
that $\al F(y)$ and $\overline{\al F(\bar y^{-1})}$, and
hence $\al F(y)$ and $\al F(y^{-1})$, are linearly independent
solutions to the eigenvalue equation $L\, \xi=z\, \xi$ for $y^2\not=1$.
Now Lemma \ref{lem:relsol} implies that $\al u(z)$ and $\al v(z)$ are
linearly independent solutions to $L\, \xi=z\, \xi$.
Since $\al u(z),\al v(z)\in S^+(z)$
we see that the deficiency indices of $L$ are $(1,1)$ in case
(1) of Lemma \ref{lem:symmform}.
For $\psi(z)$ in \eqref{eq:defGreen} we have a choice
$\psi(z)=A\al u(z) + B\al v(z)$, where we have to choose
$A,B\in\C$ such that
$\overline{A\al u(\bar z) + B\al v(\bar z)}\in {\mathcal D}_\te$.
In order to determine the possible choices
for $A,B\in\C$ we use Lemma \ref{lem:allWronskians}.

\begin{Lem}\label{lem:psiinDtheta} Let $\la_0=1-\sqrt{2}$, then
$\psi(z)=A\al u(z) + B\al v(z)\in {\mathcal D}_\te$
for
\begin{equation*}
\left( \begin{matrix} \bar A\\ \bar B \end{matrix} \right)
= \left( \begin{matrix} E & F \\ F & E\end{matrix}\right)
\left( \begin{matrix} e^{i\te} \\ e^{-i\te} \end{matrix} \right)
\qquad
\left\{ \begin{matrix}
E  = (i\la_0q^\hf e^{i\psi},-i\la_0q^\hf e^{-i\psi},
rq^{-\hf}/\la_0, q^{\frac32}\la_0/r;q)_\infty, \hfill \\
F  = (-i\la_0q^\hf e^{i\psi},i\la_0q^\hf e^{-i\psi},
-rq^{-\hf}/\la_0, -q^{\frac32}\la_0/r;q)_\infty .
\end{matrix}\right.
\end{equation*}
\end{Lem}

Of course, $A$ and $B$ are determined only up to
a common scalar. Note that $E,F\in\R$ in
Lemma \ref{lem:psiinDtheta}, and hence $\bar A=B$.
In this case we have, using \eqref{eq:conjukzreltovkz},
\begin{equation}\label{eq:psiisnotreal}
\begin{split}
\psi_k(z) &= A \al_k u_k(z) + \bar A \al_k v_k(z) =
A \al_k u_k(z) + \frac{\te(\bar t)}{\te(t)}e^{-2i\ga}
\overline{ A \al_k u_k(\bar z)} \\ &=
\frac{e^{-i\ga}|A|}{\te(t)} \Bigl(
e^{i(\ga+\arg A)}\te(t)  \al_k u_k(z)
+ \overline{e^{i(\ga+\arg A)}\te(t) \al_k u_k(\bar z)}\Bigr).
\end{split}
\end{equation}

\begin{proof} Note that $\hf(i\la_0+ (i\la_0)^{-1})=i$
and $|i\la_0|<1$, so
$\Psi(i)=e^{i\ga}\al F(i\la_0) =
\overline{e^{i\ga} \al F(-i\la_0})$ and
$\Psi(-i)= e^{i\ga}\al F(-i\la_0) =
\overline{e^{i\ga} \al F(i\la_0)}$, so we can now relate
$A$ and $B$ to the self-adjoint extension
$(L^\ast, {\mathcal D}_\te)$, see \eqref{eq:defDtheta},
by
\begin{equation}\label{eq:defsferfunkt}
\begin{split}
 0 = &\lim_{N\to\infty} [\overline{\psi(\bar z)},
e^{i\te}\Psi(i) + e^{-i\te}\Psi(-i)]_N
=  \lim_{N\to\infty} [
\overline{A\,\al u(\bar z) +B\, \al v(\bar z)},
e^{i\te}\Psi(i) + e^{-i\te}\Psi(-i)]_N \\
=& \frac{q^{\hf}e^{i\ga}}{ir}\Bigl( \bar B
\bigl\{ e^{i\te} d(i\la_0;-a,t) +
e^{-i\te} d(-i\la_0;-a,t)\bigr\}
-  \bar A\bigl\{ e^{i\te} d(i\la_0;a,t)+
e^{-i\te} d(-i\la_0;a,t)\bigr\}\Bigr)
\end{split}
\end{equation}
using Lemma \ref{lem:allWronskians}.
The condition of \eqref{eq:defsferfunkt}
determines $A$ and $B$ uniquely in terms of $e^{i\te}$ up to
a common scalar constant. Observe that all functions
$d$ in \eqref{eq:defsferfunkt}  have a common
denominator, so that we can take $A$ and $B$
as in the lemma.
\end{proof}

With $A$ and $B=\bar A$ determined by
Lemma \ref{lem:psiinDtheta}
in terms of $e^{i\te}$
we can determine the resolvent operator $R(z)$ for the
corresponding self-adjoint extension $(L^\ast, {\mathcal D}_\te)$.
For the Green kernel, see \eqref{eq:defGreen}, we need
the Wronskian, with $z=\hf(y+y^{-1})$, $|y|<1$,
\begin{equation}\label{eq:WronskianPsipsi}
[\Psi(z), \overline{\psi(\bar z)}] = e^{i\ga}
\bigl\{ \overline{A\,  c(\bar y^{-1};a,t)}
+ A\, \overline{c(\bar y^{-1};-a,t)}\bigr\}
 \hf (y^{-1}-y)
\end{equation}
by Lemmas \ref{lem:relsol} and \ref{lem:allWronskians}.

We use the parametrisation $z=\hf(y+y^{-1})$, $|y|<1$,
for $\CR$, and we want to take $z\to x\in \R$ in order
to use \eqref{eq:StieltjesPerron} to determine the spectral
measure of the self-adjoint extension $(L^\ast,{\mathcal D}_\te)$.
Note that $z\in [-1,1]$ corresponds to $y$ on the unit circle and
$z\in(-\infty,-1]$, respectively $[1,\infty)$ corresponds
to $y\in [-1,0)$, respectively in $(0,1]$. Letting $z$
tend to $x=\hf(y_0+y_0^{-1})\in\R\backslash [-1,1]$, $y_0\in(-1,1)$
from the upper or lower half plane both correspond to $y\to y_0$.
However, for $x\in(-1,1)$, put $x=\cos\chi$ with $0<\chi<\pi$,
for $\ep\downarrow0$, $z=x-i\ep \to x$ corresponds to
$y\to e^{i\chi}$ and $z=x+i\ep \to x$ corresponds to
$y\to e^{-i\chi}$. So we consider these cases separately.
For the moment we assume $\xi,\eta\in{\mathcal D}$ so that all
summations are actually finite, the general case
$\xi,\eta\in\lt$ follows by continuity of the spectral
projections $E({\mathcal B})$, ${\mathcal B}\subset \R$
a Borel set.

\begin{Prp}\label{prop:continuouspartL} $[-1,1]$ is contained
in the continuous spectrum of $(L,{\mathcal D}_\te)$ and
for $0\leq\chi_1<\chi_2\leq\pi$ the spectral measure
is determined by
\begin{equation*}
\langle E\bigl( [\cos\chi_2,\cos\chi_1]\bigr) \xi, \eta \rangle =
\frac{1}{2\pi} \int_{\chi_1}^{\chi_2}
\frac{\langle \xi, \psi(\cos\chi)\rangle
\overline{\langle \eta, \psi(\cos\chi) \rangle} }
{ \bigl| A\, c(e^{i\chi};a,t)
+ \bar A \, c(e^{i\chi};-a,t)\bigr|^2}\, d\chi
\end{equation*}
where $A$, $B=\bar A$ is determined by
Lemma \ref{lem:psiinDtheta}.
\end{Prp}

\begin{proof}
We first assume $x\in(-1,1)$, $x=\cos\chi$, $0<\chi<\pi$.
Observe that for $\ep\downarrow 0$, $\Psi_k(x-i\ep)\to
e^{i\ga}\al_k F_k(e^{i\chi})$, $\Psi_k(x+i\ep)\to
e^{i\ga}\al_k F_k(e^{-i\chi})$, and
$\psi_k(x\pm i\ep)\to  \psi_k(\cos\chi)$,
so by \eqref{eq:exprresolvinGreen}
and Lemma \ref{lem:allWronskians}
\begin{equation}\label{eq:limGreen}
 \lim_{\ep\downarrow 0}
\langle R(x-i\ep)\xi,\eta\rangle - \langle R(x+i\ep)\xi,\eta\rangle =
2\sum_{k\leq l} A_{k}
\frac{\overline{\psi_l(\cos\chi)}}{(e^{-i\chi}-e^{i\chi})}
\bigl( \xi_l\bar \eta_k+\xi_k\bar \eta_l\bigr) \bigl(1-\hf
\de_{k,l}\bigr)
\end{equation}
with
\begin{equation*}
\begin{split}
A_{k} = &
\frac{\al_k F_k(e^{i\chi})}
{\bigl( \overline{A\, c(e^{i\chi};a,t)}
+ A\, \overline{c(e^{i\chi};-a,t)}\bigr)}
+ \frac{\al_k F_k(e^{-i\chi})}
{\bigl(\overline{ A\, c(e^{-i\chi};a,t)}
+ A\, \overline{c(e^{-i\chi};-a,t)}\bigr)} \\
= &  \frac{\bigl(\overline{ A\, c(e^{-i\chi};a,t)}
+ A\, \overline{c(e^{-i\chi};-a,t)}\bigr)\al_k F_k(e^{i\chi})
+ \bigl( \overline{A\, c(e^{i\chi};a,t)}
+ A\, \overline{c(e^{i\chi};-a,t)}\bigr)\al_k F_k(e^{-i\chi})}
{\bigl( \overline{A\, c(e^{i\chi};a,t)}
+ A\, \overline{c(e^{i\chi};-a,t)}\bigr)
\bigl(\overline{ A\, c(e^{-i\chi};a,t)}
+ A\, \overline{c(e^{-i\chi};-a,t)}\bigr)} \\
= & \frac{\bigl(\bar A\, \al_k v_k(\cos\chi) + A \, \al_k
u_k(\cos\chi)\bigr) }
{\overline{\bigl( A\, c(e^{i\chi};a,t)
+ \bar A \, c(e^{i\chi};-a,t)\bigr)}
\bigl(\bar A\, c(e^{i\chi};-a,t)
+ A \, c(e^{i\chi};a,t)\bigr)} \\
= & \frac{\psi_k(\cos\chi)}
{\Bigl| A\, c(e^{i\chi};a,t)
+ \bar A \, c(e^{i\chi};-a,t)\bigr|^2}
\end{split}
\end{equation*}
using \eqref{eq:relforcf} and Lemma \ref{lem:relsol}
for the third equality.
The above gives an explicit expression for $A_k$.
If the expression $A_k\overline{\psi_l(\cos\chi)}$
is symmetric in $k$ and $l$ we can rewrite the sum
over $k\leq l$ in \eqref{eq:limGreen} as the product of a sum
over $k$ and a sum over $l$. From
\eqref{eq:psiisnotreal} we see that
$\psi_k(x)\overline{\psi_l(x)}$ is symmetric in $k$ and $l$
for $x\in\R$, and hence $A_k\overline{\psi_l(\cos\chi)}$
is symmetric in $k$ and $l$.

So we can antisymmetrise the sum in \eqref{eq:limGreen}
and  \eqref{eq:limGreen} equals
\begin{equation*}
\frac{\Bigl( \sum_{l=-\infty}^\infty \xi_l \overline{\psi_l(\cos\chi)}
\Bigr)
\Bigl( \sum_{k=-\infty}^\infty \psi_k(\cos\chi) \overline{\eta_k}
\Bigr)}
{\bigl| A\, c(e^{i\chi};a,t)
+ \bar A \, c(e^{i\chi};-a,t)\bigr|^2}.
\end{equation*}
Using $dx = \frac{1}{2i}(e^{i\chi}-e^{-i\chi})d\chi$,
dominated convergence and
\eqref{eq:StieltjesPerron} we find for $0<\chi_1<\chi_2<\pi$
\begin{equation*}
\langle E\bigl( (\cos\chi_2),\cos(\chi_1)\bigr) \xi, \eta \rangle =
\frac{1}{2\pi} \int_{\chi_1}^{\chi_2}
\frac{\langle \xi, \psi(\cos\chi)\rangle
\overline{\langle \eta, \psi(\cos\chi) \rangle} }
{ \bigl| A\, c(e^{i\chi};a,t)
+ \bar A \, c(e^{i\chi};-a,t)\bigr|^2}\, d\chi.
\end{equation*}

By the previous calculation the proposition
follows for the open interval $(-1,1)$. Since the spectrum
is closed we see that $\pm 1$ are contained in the spectrum
of $(L^\ast,{\mathcal D}_\te)$. Since $(L^\ast,{\mathcal D}_\te)$
is a self-adjoint operator $\pm 1$ can be in the continuous
spectrum or in the point spectrum, see \cite[Thm.~13.27]{DunfS}.
The proposition follows by showing that the endpoints $\pm 1$
are not contained in the discrete spectrum.
Note that for $y=\pm 1$ the first Wronskian
in Lemma \ref{lem:allWronskians} vanishes, so we need to
construct a second independent
solution to $\al F(\pm 1)$ first. We consider the
case $y=1$, the case $y=-1$ is being dealt with
similarly. Put
\begin{equation}\label{eq:defothersol}
H_k(y) = \frac{F_k(y)-F_k(1)}{y-1},
\end{equation}
then it satisfies
$$
\frac{1+a^2tq^{k-1}}{atq^{k-1}}\, H_{k+1}(y)
- \frac{1-q^{k-1}t}{atq^{k-1}}\, H_{k-1}(y)
= 2\, H_k(y) + (1-y^{-1})\, F_k(y).
$$
Taking $y\to 1$ gives the solution $H_k(1) =
\frac{\partial F_k(y)}{\partial y}\big\vert_{y=1}$
to \eqref{eq:diffequnnorm} for the eigenvalue $z=1$.
Now \eqref{eq:defothersol} gives the asymptotic
behaviour
\begin{equation}\label{eq:dasymothersol}
H_k(1) = \frac{\partial F_k(y)}{\partial y}\big\vert_{y=1} =
(-k)a^{-k}\bigl( 1+ {\mathcal O}(q^{-k})\bigr),
\qquad k\to-\infty.
\end{equation}
Using the asymptotic behaviour \eqref{eq:dasymothersol}
we can calculate the Wronskian
\begin{equation*}
[\al H(1),  \overline{\al F(1)}] = -\hf
\end{equation*}
similar to the calculation of the first
Wronskian of Lemma \ref{lem:allWronskians}.
So we have two linearly independent solutions of the
eigenvalue equation $L^\ast\, \xi= \xi$. From the asymptotic
behaviour \eqref{eq:dasymothersol} of $\al H_k$ and
of $\overline{\al F(1)}$ as $k\to-\infty$, it follows
that no linear combination of $\al H$ and $\overline{\al F(1)}$
can be an element of $\lt$. Hence, $1$ is not in the point
spectrum, and hence $1$ is contained in the continuous spectrum.
\end{proof}

Note that with the choices for $\Psi(z)$ and $\psi(z)$
the function $[\Psi(z), \overline{\psi(\bar z)}]\, G_{k,l}(z)$
is analytic for $z\in\C\backslash[-1,1]$, and the same
holds for
$[\Psi(z), \overline{\psi(\bar z)}]\, \langle R(z)\xi,\eta\rangle$,
where we still assume $\xi,\eta\in{\mathcal D}$.
We next turn to the spectrum of $(L,{\mathcal D}_\te)$
contained in $(-\infty,-1)\cup(1,\infty)$.
Because of these remarks and the remarks in
the paragraph preceding
Proposition \ref{prop:continuouspartL} and
\eqref{eq:exprresolvinGreen} we see that
$E_{\xi,\eta}\bigl( (x_1,x_2)\bigr)=0$ as long as
$(x_1,x_2)$ contains no zero of the
Wronskian  \eqref{eq:WronskianPsipsi}
using dominated convergence in \eqref{eq:StieltjesPerron}.
Note that the zeroes of the Wronskian are isolated, since
the Wronskian \eqref{eq:WronskianPsipsi} is
meromorphic in $y$. So the only discrete mass points can
occur at a zero of the Wronskian  \eqref{eq:WronskianPsipsi}.

\begin{Prp}\label{prop:discretespectrum}
There is no continuous spectrum
of $(L^\ast, {\mathcal D}_\te)$ in $(-\infty,-1)\cup(1,\infty)$.
The point spectrum of $(L^\ast, {\mathcal D}_\te)$ occurs at
the set
\begin{equation*}
S = \Bigl\{ x_0 =\hf(y_0+y_0^{-1}) \,\Bigl\vert\,
|y_0|>1, \  \bar A c(y_0;-a,t) + A c(y_0;a,t)=0\Bigr\}
\end{equation*}
and the spectral projection is determined by
\begin{equation*}
\langle E(\{x_0\})\xi,\eta\rangle =
\text{{\rm Res}}\, 
\frac{\langle \xi, \psi(x_0)\rangle \langle \psi(x_0), \eta\rangle}
{y\bigl(A c(y^{-1};a,t)+\bar A c(y^{-1};-a,t)\bigr)
\bigl(\bar A c(y;-a,t)+  A c( y;a,t)\bigr)}\Bigl\vert_{y=y_0}.
\end{equation*}
\end{Prp}

\begin{Rmk}\label{rmk:propdiscretespectrum}
From Proposition \ref{prop:discretespectrum} and the
discussion preceding it we see
that the discrete set $S$ is contained in $\R$.
Moreover, the unboudedness of $(L^\ast,{\mathcal D}_\te)$
and the boundedness of the continuous spectrum $[-1,1]$,
see Proposition \ref{prop:continuouspartL}, shows that the set $S$
is unbounded.
\end{Rmk}

\begin{proof} From the remarks preceding Proposition
\ref{prop:discretespectrum} we see that we can only
have discrete spectrum in $(-\infty,-1)\cup(1,\infty)$.
Next assume that $x_0\in (-\infty,-1)\cup(1,\infty)$ is
a zero of the Wronskian  \eqref{eq:WronskianPsipsi}.
(Note that we have already dealt with the case $x_0=\pm 1$
in Proposition \ref{prop:continuouspartL}.)
Let $x_0=\hf(y_0+y_0^{-1})$ with $|y_0|>1$. (Note that
this is against the convention, but it makes formulas
better looking.)
Moreover, since $\langle R(z)\xi,\eta\rangle$ is meromorphic
in a neighbourhood of $x_0$ we find that $x_0$ is an
element of the point spectrum of $(L, {\mathcal D}_\te)$
and
\begin{equation}\label{eq:discrspectrum1}
\langle E(\{x_0\})\xi,\eta\rangle = \frac{1}{2\pi i}
\oint_{\mathcal C} \langle R(z)\xi,\eta\rangle \, dz =
\text{{\rm Res}}\, \langle R(z)\xi,\eta\rangle\bigl\vert_{z=x_0}
\end{equation}
where $\mathcal C$ is a small positively oriented
contour enclosing $x_0$ once and no other
singularities of the resolvent. From
Lemma \ref{lem:relsol} and the fact that $y_0^{-1}$  is a
zero of the Wronskian \eqref{eq:WronskianPsipsi} using
\eqref{eq:relforcf} we find
\begin{equation*}
\begin{split}
\psi_k(x_0) &= A\al_k u_k(x_0) + \bar A\al_k v_k(x_0)
= \bigl( A c(y_0^{-1};a,t) +\bar A c(y_0^{-1};-a,t) \bigr) \, \al_k F_k(y_0^{-1}) \\
&= \bigl( A c(y_0^{-1};a,t) +\bar A c(y_0^{-1};-a,t) \bigr) \,
e^{-i\ga} \Psi_k(x_0).
\end{split}
\end{equation*}
In particular, this implies $\psi(x_0)\in\lt$.
Using this in \eqref{eq:discrspectrum1} and switching
from $z$ to $y$ gives the desired expression for
$\langle E(\{x_0\})\xi,\eta\rangle$.
\end{proof}

Combining Propositions \ref{prop:continuouspartL} and
\ref{prop:discretespectrum} proves the following
theorem, which summarises the results of this
section.

\begin{Thm} \label{thm:case1spectraldecomp}
The spectral decomposition of the self-adjoint
extension $(L^\ast, {\mathcal D}_\te)$
defined by \eqref{eq:defDtheta} of $(L,{\mathcal D})$
as defined in \eqref{eq:defLonlt} with $a_k$
as in \eqref{eq:ak1} and $a=q^\hf e^{i\psi}$,
$t=ire^{-i\psi}\notin\R_{>0}$, $r\in\R\backslash\{0\}$, is
given by
\begin{equation*}
\begin{split}
\langle L^\ast \xi,\eta\rangle =&
\frac{1}{2\pi} \int_0^\pi
\frac{\cos\chi\, ({\mathcal F_\te}\xi)(\cos\chi) \overline{({\mathcal
F_\te}\eta)(\cos\chi)}}
{ \bigl| A\, c(e^{i\chi};a,t)
+ \bar A \, c(e^{i\chi};-a,t)\bigr|^2}\, d\chi
\\ &+ \sum_{x_0\in S}
\text{{\rm Res}}\, 
\frac{x_0\, ({\mathcal F_\te}\xi)(x_0) \overline{({\mathcal
F_\te}\eta)(x_0)}}
{y\bigl(A c(y^{-1};a,t)+\bar A c(y^{-1};-a,t)\bigr)
\bigl(\bar A c(y;-a,t)+  A c( y;a,t)\bigr)}\Bigl\vert_{y=y_0}
\end{split}
\end{equation*}
where $\xi\in{\mathcal D}_\te$, $\eta\in\lt$,
$A$, $B=\bar A$ is determined in
Lemma \ref{lem:psiinDtheta} by $e^{i\te}$,
$c(\cdot;a,t)$ is defined in Lemma \ref{lem:relsol}, the
set of discrete mass points is given by
\begin{equation*}
S = \Bigl\{ x_0 =\hf(y_0+y_0^{-1}) \,\Bigl\vert\,
|y_0|>1, \  \bar A c(y_0;-a,t) + A c(y_0;a,t)=0\Bigr\},
\end{equation*}
the corresponding Fourier transform is
\begin{equation*}
({\mathcal F_\te}\xi)(x) = \langle \xi, \psi(x)\rangle =
\langle \xi, A\al u(x) + \bar A\al v(x)\rangle
\end{equation*}
with $\al_k$ defined by \eqref{eq:ak1}, and
$u_k(x), v_k(x)$ as in Lemma \ref{lem:relsol}.
\end{Thm}

We recast Theorem \ref{thm:case1spectraldecomp} into
two immediate corollaries.

\begin{Cor}\label{cor:thmcase1spectraldecomp}
With the notation of
Theorem \ref{thm:case1spectraldecomp}
the orthogonality relations
\begin{equation*}
\begin{split}
&\de_{k,l} = \frac{1}{2\pi} \int_0^\pi
\frac{\psi_k(\cos\chi) \overline{\psi_l(\cos\chi)}}
{ \bigl| A\, c(e^{i\chi};a,t)
+ \bar A \, c(e^{i\chi};-a,t)\bigr|^2}\, d\chi
\\ &+ \sum_{x_0\in S}
\text{{\rm Res}}\, 
\frac{ \psi_k(x_0) \overline{\psi_l(x_0)}}
{y\bigl(A c(y^{-1};a,t)+\bar A c(y^{-1};-a,t)\bigr)
\bigl(\bar A c(y;-a,t)+  A c( y;a,t)\bigr)}\Bigl\vert_{y=y_0}
\end{split}
\end{equation*}
hold, and the functions $\bigl\{\psi_k\bigr\}_{k\in\Z}$
form an orthonormal basis of the corresponding
weighted $L^2$-space, and ${\mathcal F}_\te$ is
a unitary isomorphism from $\lt$ to the
corresponding weighted $L^2$-space.
\end{Cor}

\begin{Cor}\label{cor2:thmcase1spectraldecomp} With the notation of
Theorem \ref{thm:case1spectraldecomp} we have
the following transform pair; for
$\xi=\sum_{k=-\infty}^\infty \xi_ke_k\in\lt$
\begin{equation*}
\begin{split}
&\xi_l = \frac{1}{2\pi} \int_0^\pi
\frac{({\mathcal F}_\te \xi)(\cos\chi) \overline{\psi_l(\cos\chi)}}
{ \bigl| A\, c(e^{i\chi};a,t)
+ \bar A \, c(e^{i\chi};-a,t)\bigr|^2}\, d\chi
\\ &+ \sum_{x_0\in S}
\text{{\rm Res}}\, 
\frac{  ({\mathcal F}_\te \xi)(x_0)\overline{\psi_l(x_0)}}
{y\bigl(A c(y^{-1};a,t)+\bar A c(y^{-1};-a,t)\bigr)
\bigl(\bar A c(y;-a,t)+  A c( y;a,t)\bigr)}\Bigl\vert_{y=y_0}.
\end{split}
\end{equation*}
\end{Cor}

\begin{Rmk}\label{rmk:discspectralE}
Note that the spectral measure for the continuous
spectrum in Theorem \ref{thm:case1spectraldecomp} is rather explicit,
and that the description
of the discrete mass points in Theorem \ref{thm:case1spectraldecomp}
is indirect. For
the special case of the Ismail-Zhang
$q$-analogue of the exponential function defined
in \eqref{eq:defcurlyE} we can describe the discrete mass
points a bit more explicitly.
This special case corresponds to $\psi=0 \mod \pi$.
Without loss of generality we can assume $\psi=0$ by
\cite[(1.4.6)]{GaspR}, and take $a=q^\hf$ and
$t=ir$, $r\in\R\backslash\{0\}$. First observe that in
the definition of $A=e^{i\te}E+e^{-i\te}F$
in Lemma \ref{lem:psiinDtheta} we can replace
$E$ and $F$ by, recall $\la_0=1-\sqrt{2}$,
$E  = \te(rq^{-\hf}/\la_0)$ and
$F  = \te(-rq^{-\hf}/\la_0)$ by cancelling a
common factor.
In this case the $c$-functions
have a common factor and we have
\begin{equation}\label{eq:weightE1}
A\, c(y;q^\hf,ir) +\bar A \, c(y;-q^\hf,ir) =
\frac{(q^\hf/y,-q^\hf/y;q)_\infty}{(-q;q)_\infty \te(ir)(y^{-2};q)_\infty}
\Bigl( \te(yq^\hf ir)A
+ \te(-yq^\hf ir)\bar A\Bigr).
\end{equation}
So the spectral measure for the continuous part can be read off
from \eqref{eq:weightE1}. For the discrete spectrum we have
to find the zeroes of \eqref{eq:weightE1} as function of $y$
for $|y|>1$, so we have to solve
$\te(yq^\hf ir)A  = - \te(-yq^\hf ir)\bar A$. Put $y=e^{2\pi iw}$
and consider
\begin{equation}\label{eq:ellipticfunctiong}
g(w,\tau) = \frac{\te(e^{2\pi iw}q^\hf ir)}
{\te(-e^{2\pi iw}q^\hf ir)}, \qquad q=e^{\pi i\tau},\quad \tau\in
i\R_{>0},
\end{equation}
so that the equation is rewritten as $g(w,\tau)= -\bar A/A$.
It follows from \eqref{eq:thetaidentity} that
$g(w,\tau)$ is an elliptic function with periods $1$ and
$\tau$. From \cite[Ch. XX, XXI]{WhitW} we see that the
order of the elliptic function $g$ is $2$, so that the
equation $g(w,\tau)= -\bar A/A$ has $2$ solutions in each
fundamental parallellogram. Since the solutions in the
$y$-coordinate are
real we find $w\in i\R\cup \hf+i\R$ (modulo $1$).
By period $\tau$
it follows that the discrete mass points are of the
form $x^{(i)}_n = \hf(y_iq^{-2n}+y_i^{-1}q^{2n})$, $n\in\Z$,
with $|x^{(i)}_n|>1$, where $i=1,2$. So in particular, in this
case the discrete mass points are located
on two $q^2$-quadratic grids.

In this case we have from \eqref{eq:psiisnotreal}
\begin{equation}
\begin{split}
\psi_k(x) &=
\frac{2e^{-i\ga}|A|}{\te(ir)}
\Re \bigl[ q^{\hf k} e^{i(\phi_k+\ga+\arg A)}
\frac{\te(ir) (-r^2q^{2+4k};q^4)_\infty}{(irq^k;q)_\infty}
\E_{q^2}(x;-irq^k)\bigr], \qquad x\in\R,
\end{split}
\end{equation}
so that we can look upon the orthogonality relations
of Corollary \ref{cor:thmcase1spectraldecomp}
or the integral transform of
Corollary \ref{cor2:thmcase1spectraldecomp}
as a $q$-analogue
of the Fourier cosine transform
for the $q$-exponential $\E_q$ as defined in
\eqref{eq:defcurlyE} using \eqref{eq:limcurlyE}.
\end{Rmk}

\begin{Rmk} \label{rmk:explicittransforms}
It is of interest to be able to calculate the
${\mathcal F}_\te$ transforms of specific vectors
and next use Corollary \ref{cor2:thmcase1spectraldecomp}
to get explicit
transforms, even though the precise location of the
discrete mass points is not known.
Results already present in the literature can be
used for this. Since the corresponding formulas are
well known we leave it to the reader to fill in the
details. As a first example, the
${\mathcal F}_\te$-transform of the vector
$\xi = \sum_{k=-\infty}^\infty z^k \al_k\, e_k$ can
be expressed in terms of
infinite $q$-shifted factorials
using the generating function
\cite[Lemma~3.3 with
$k=1$]{KoelR}, \cite[Lemma~2.2 with
$k=1$]{Rose}. Since this is not an $\lt$-vector some care
has to be taken, but using an approximation
argument plus the absolute convergence of the sum
defining ${\mathcal F}_\te \xi$
for $z$ in a certain annulus, we can find the result.

Using a generalisation of Rahman's summation
formulas, see \cite[Prop.~3.1]{KoelR},
\cite[Thm.~2.1 with $k=l=1$]{Rose} it is possible
to calculate the Poisson kernel, i.e. the
 ${\mathcal F}_\te$-transform of the vector
$\sum_{k=-\infty}^\infty z^k \psi_k(x^\prime)\, e_k$
for a different value for the argument.
For the arguments in the interval $[-1,1]$
the Poisson kernel can be expressed in term
of eight very-well-poised ${}_8W_7$-series.
For the case $a=q^\hf$, $t=ir$,
i.e. for the situation corresponding to
$\E_{q^2}$, the situation
simplifies greatly, and the eight ${}_8W_7$-series
can be combined to only two ${}_8W_7$-series
by \cite[(2.10.1)]{GaspR}.
To evaluate the
Poisson kernels in the discrete mass points
we express $\psi_k(x_0)$ in terms of $\al F(y_0^{-1})$
as a single ${}_2\vp_1$,
and use the connection coefficients of Lemma \ref{lem:relsol}
before applying the same summation formulas again.
The procedure sketched above can be generalised to
the ${\mathcal F}_\te$-transform of the vector
$$
\sum_{k=-\infty}^\infty z^k
\rvps{r+1}{r}{a_1,\ldots,a_{r+1}}{b_1,\ldots,b_r}{tq^k} \, e_k
$$
using \cite[Thm.~2.1, with $k=r$, $l=1$]{Rose}.
\end{Rmk}


\section{Spectral decomposition of $L$ in the second case}
\label{sec:specdecL2}

In this section we calculate the spectral measure as
explicitly as possible of the self-adjoint extensions of
$(L,{\mathcal D})$ with $L$ as in \eqref{eq:defLonlt} with $a_k$
given by \eqref{eq:ak2}. As in \S \ref{sec:specdecL2},
this depends on the parameter $\te$
of the self-adjoint extension $(L^\ast,{\mathcal D}_\te)$ of
$(L,{\mathcal D})$. There is always continuous spectrum on the
interval $[-1,1]$, and an infinite series of discrete mass
points tending to plus or minus $\infty$. The location of the
discrete mass points depends on the choice of the
self-adjoint extension. In this section we always
have $a=is$, $s\in\R\backslash\{0\}$ and $t<0$ as
in case (2) of Lemma \ref{lem:symmform}, but we
keep the notation $a$ and $t$ in order to keep the analogy with
\S \ref{sec:specdecL1}.
The case considered in this section is slightly easier than
the case considered in \S \ref{sec:specdecL1}, so
we only state the results and indicate the proofs by analogy
to \S \ref{sec:specdecL1}.

So, in this section $t<0$, $a=is$, $s\in\R\backslash\{0\}$,
see case (2) of Lemma \ref{lem:symmform} and $a_k$ and $\al_k$
are given in \eqref{eq:ak2}.
In this case it is straightforward to
see that $\overline{\al_kF_k(\bar y)}=\al_k F_k(y)$.
Using Lemma \ref{lem:relsol} and \eqref{eq:relforcf} this
implies
\begin{equation}\label{eq:baruisv}
\overline{\al_k\,  u_k (\bar z)} = \al_k\, v_k(z).
\end{equation}
We also have that $\al_kF_k(y) = y^{-k}(1+{\mathcal O}(q^{-k}))$,
so that $S^-(z)$ is spanned by $\Psi(z)=\al F(y)$ with
$z=\hf(y+y^{-1})$ with $|y|<1$. The
statement analogous to Lemma \ref{lem:allWronskians} is
the following lemma, whose proof is similar to the proof of
Lemma \ref{lem:allWronskians}.

\begin{Lem}\label{lem:allWronskians2} We have the following Wronskians;
\begin{equation*}
\begin{split}
[\al F(y), \al F(y^{-1})] &= \hf (y^{-1}-y), \\
\lim_{N\to\infty}
[\overline{\al u(\bar w)}, \al F(y)]_N  &=
\frac{i q}{st} \frac{\te(s^2t/q)}{\te(t)}
\,  d(y;a,t), \\
\lim_{N\to\infty}
[\overline{\al v(\bar w)}, \al F(y)]_N &=
-\frac{i q}{st} \frac{\te(s^2t/q)}{\te(t)}
\,  d(y;-a,t).
\end{split}
\end{equation*}
\end{Lem}

It follows from Lemma \ref{lem:allWronskians2}
that $\al F(y)$ and $\overline{\al F(\bar y^{-1})}$, and
hence $\al F(y)$ and $\al F(y^{-1})$, are linearly independent
solutions to the eigenvalue equation $L\, \xi=z\, \xi$ for $y^2\not=1$.
Now Lemma \ref{lem:relsol} implies that $\al u(z)$ and $\al v(z)$ are
linearly independent solutions to $L\, \xi=z\, \xi$.
Since $\al u(z),\al v(z)\in S^+(z)$
we see that the deficiency indices of $L$ are $(1,1)$ in case
(2) of Lemma \ref{lem:symmform}.
Again, for $\psi(z)$ in \eqref{eq:defGreen} we have a choice
$\psi(z)=A\al u(z) + B\al v(z)$, where we have to choose
$A,B\in\C$ such that
$\overline{A\al u(\bar z) + B\al v(\bar z)}\in {\mathcal D}_\te$.
In order to determine the possible choices
for $A,B\in\C$ we use Lemma \ref{lem:allWronskians2}.

\begin{Lem}\label{lem:psiinDtheta2} Let $\la_0=1-\sqrt{2}$, then
$\psi(z)=A\al u(z) + B\al v(z)\in {\mathcal D}_\te$
for
\begin{equation*}
\left( \begin{matrix} \bar A\\ \bar B \end{matrix} \right)
= \left( \begin{matrix} E & F \\ F & E\end{matrix}\right)
\left( \begin{matrix} e^{i\te} \\ e^{-i\te} \end{matrix} \right)
\qquad
\left\{ \begin{matrix}
E  = (s\la_0,\la_0q/s,-st/q\la_0,-q^2\la_0/st;q)_\infty, \hfill \\
F  = (-s\la_0,-\la_0q/s,st/q\la_0,q^2\la_0/st;q)_\infty.
\end{matrix}\right.
\end{equation*}
\end{Lem}

Of course, $A$ and $B$ are determined only up to
a common scalar. Note that $E,F\in\R$ in
Lemma \ref{lem:psiinDtheta2}, and hence $\bar A=B$.
In this case we have, using \eqref{eq:baruisv},
\begin{equation}\label{eq:psiisreal2}
\begin{split}
\psi_k(z) &= A \al_k u_k(z) + \bar A \al_k v_k(z) =
A \al_k u_k(z) + \overline{ A \al_k u_k(\bar z)},
\end{split}
\end{equation}
so that for $z=x\in\R$ we have
$\psi_k(x)=2\Re [A \al_k u_k(x)]$ is real-valued.

Completely analagous to Proposition
\ref{prop:continuouspartL} we obtain that
$[-1,1]$ is  contained
in the continuous spectrum of $(L,{\mathcal D}_\te)$ and
for $0\leq\chi_1<\chi_2\leq\pi$ the spectral measure
is determined by the same formula as in
Proposition
\ref{prop:continuouspartL}, but with $a=is$, $t<0$ and
where $A$, $B=\bar A$ are determined by
Lemma \ref{lem:psiinDtheta2}.

The expression \eqref{eq:WronskianPsipsi} for
the Wronskian has to be replaced by, again
$z=\hf(y+y^{-1})$, $|y|<1$,
\begin{equation}\label{eq:Wronskaincase2}
[\Psi(z), \overline{\psi(\bar z)}] =
\bigl\{ \overline{A\,  c(\bar y^{-1};a,t)}
+ A\, \overline{c(\bar y^{-1};-a,t)}\bigr\}
 \hf (y^{-1}-y),
\end{equation}
and since for the discrete spectrum only the zeroes
of the Wronskian play a role, we see that
Proposition \ref{prop:discretespectrum}
goes through in this case
with $A$ and $\bar A=B$ defined by Lemma
\ref{lem:psiinDtheta2} in this case.

Combining these results then gives the spectral
decomposition of the self-adjoint extension
$(L^\ast,{\mathcal D}_\te)$ of $(L,{\mathcal D})$
as in \eqref{eq:defLonlt} with $a_k$ defined
by \eqref{eq:ak2}.

\begin{Thm} \label{thm:case2spectraldecomp}
The spectral decomposition of the self-adjoint
extension $(L^\ast, {\mathcal D}_\te)$
defined by \eqref{eq:defDtheta} of $(L,{\mathcal D})$
as defined in \eqref{eq:defLonlt} with $a_k$
as in \eqref{eq:ak1} and $a=is$, $s\in\R\backslash\{0\}$,
$t\in\R_{<0}$, is given by the same formula as in
Theorem \ref{thm:case1spectraldecomp} except that
$A$ and $\bar A=B$ are defined by Lemma
\ref{lem:psiinDtheta2}.
Moreover, Corollary \ref{cor:thmcase1spectraldecomp}
and Corollary \ref{cor2:thmcase1spectraldecomp}
remain valid in this case.
\end{Thm}


\section{Quadratic transformation}
\label{sec:quadrtrans}

In this section we relate some of the solutions
discussed in Lemma \ref{eq:solfrombigqJacobi}
to ${}_3\vp_2$-series of base $q^2$. The resulting
transformation of Proposition \ref{prop:quadrtransform}
can be considered
as a non-terminating analogue of Singh's quadratic
transformation \cite[(3.10.13)]{GaspR}. The reason
for this is that the symmetric
Al-Salam--Chihara polynomials of even degree can be expressed
in terms of continuous dual $q^2$-Hahn polynomials,
and the resulting transform is a special case of
Singh's transformation. It does not seem possible to
obtain the result of Proposition \ref{prop:quadrtransform} as
a special or limiting case of Singh's transformation.

The recurrence relation \eqref{eq:diffequnnorm}
has no term involving $f_k(z)$ in the right
hand side, so we can iterate the recurrence to
obtain a three term recurrence for the
even and odd degree $f_k(z)$'s.
For convenience put
\begin{equation*}
c_k = \frac{1+a^2t q^{k-1}}{atq^{k-1}} = a(1+q^{1-k}/at), \qquad
d_k = - \frac{1-q^{k-1}t}{atq^{k-1}} =  \frac{1}{a}
(1-q^{1-k}/t),
\end{equation*}
and iterate \eqref{eq:diffequnnorm} to find
\begin{equation}\label{eq:quadrrecrel}
(2z)^2\, f_k(z) = c_kc_{k+1}\, f_{k+2}(z) + \bigl( c_k d_{k+1} +
d_kc_{k-1}\bigr)\, f_k(z) + d_k d_{k-1}\, f_{k-2}(z).
\end{equation}

So from \eqref{eq:quadrrecrel} we find a
three-term recurrence relation for the even
and odd degree $f_k$'s.
The recurrence \eqref{eq:quadrrecrel} can be matched to
the one studied by Gupta, Ismail and Masson
\cite{GuptIM}, where a lot of solutions are discussed.
The recurrence relation \eqref{eq:quadrrecrel} is also
studied in detail in \cite{KoelSCA} as a linear
operator on a suitable Hilbert space.

We recall from \cite[\S 2-4]{KoelSCA} that
\begin{equation}\label{eq:solfrombigqJacobi}
\Phi_\ga(xq^k;a,b,c;q) =
\frac{(-q^{1-k}/bcx,-q^{1-k}\ga/ax;q)_\infty}
{(-q^{1-k}/abx,-q^{1-k}/acx;q)_\infty} (a\ga)^{-k}
\rvps{3}{2}{q\ga/a,b\ga,c\ga}{-q^{1-k}\ga/ax,q\ga^2}
{\frac{-q^{1-k}}{bcx}}
\end{equation}
is a solution of
\begin{equation}\label{eq:recrelfrombigqJacobi}
\begin{split}
&(\ga+\ga^{-1})\, F(xq^k) =
a\bigl(1+\frac{q^{-k}}{abx}\bigr)\bigl(1+\frac{q^{-k}}{acx}\bigr)\,
F(xq^{k+1})\\
&- \Bigl( q^{-k}\bigl(\frac{1}{bx}+\frac{1}{cx}+\frac{q}{abcx}+\frac{1}{ax}\bigr)
+\frac{q^{-2k}}{x^2abc}(1+q)\Bigr)\, F(xq^k)
+ \frac{1}{a} \bigl(
1+\frac{q^{1-k}}{bcx}\bigr)\bigl(1+\frac{q^{-k}}{x}\bigr)\,
F(q^{k-1}x).
\end{split}
\end{equation}
Many other solutions and their connections are known,
see \cite{GuptIM}, \cite{KoelSCA}, but we only need this
solution. The main result of this section is the
following proposition.

\begin{Prp}\label{prop:quadrtransform}
For $|z|<\min (1,|a|^2)$
\begin{equation*}\label{eq:quadrel2}
\tvpo{ay,-ay}{qy^2}{-\frac{z}{a^2}} =
\frac{(z,qzy^2/a;q^2)_\infty}
{(-z/a;q)_\infty}
\,_3\vp_2\left(
\genfrac{.}{.}{0pt}{}{q^2y^2/a,-y^2,-qy^2}{qzy^2/a,q^2y^4}
;q^2,z\right).
\end{equation*}
\end{Prp}

\begin{proof} A straightforward calculation
starting from \eqref{eq:quadrrecrel} shows
that $a^kR_k$ with $R_k=f_{2k}(z)$, $z=\hf(y+y^{-1})$,
satisfies
\eqref{eq:recrelfrombigqJacobi} in base $q^2$ with
$(a,b,c,x)$ of \eqref{eq:recrelfrombigqJacobi} specialised
to $(a,-1,-q,-t/q)$ and
$\ga=y^2$.
So the solutions of \eqref{eq:recrelfrombigqJacobi}
are related to the ones in Lemma \ref{lem:diffeqsol} for $k$
replaced by $2k$. Since the solution space is two-dimensional,
we cannot immediately give direct relations. However, the
space of subdominant or minimal solutions is one-dimensional,
see \cite[Thm.~1]{GuptIM}, and spanned by $F_{2k}(y)$, $|y|<1$,
and $a^{-k}\Phi_{y^2}(-tq^{2k-1};a,-1,-q;q^2)$. For the
renormalised recurrence in case of \S\S \ref{sec:specdecL1},
\ref{sec:specdecL2}, this is just the statement that $S^-(z)$ is
one-dimensional.

So these two
solutions only differ by a constant $C$ which can be determined
by considering the limit behaviour for $k\to-\infty$.
The limit behaviour follows from the explicit
expression in \eqref{eq:solfrombigqJacobi}
and Lemma \ref{lem:diffeqsol}.
This gives
$C=1$, and canceling common factors gives, for $k\in\Z$
and $|y|<1$, the relation
\begin{equation}\label{eq:quadrel1}
\tvpo{ay,-ay}{qy^2}{-\frac{q^{2-2k}}{a^2t}} =
\frac{(q^{2-2k}/t,q^{3-2k}y^2/at;q^2)_\infty}
{(-q^{3-2k}/at,-q^{2-2k}/at;q^2)_\infty}
\,_3\vp_2\left(
\genfrac{.}{.}{0pt}{}{q^2y^2/a,-y^2,-qy^2}{q^{3-2k}y^2/at,q^2y^4}
;q^2,\frac{q^{2-2k}}{t}\right).
\end{equation}
Multiplying
\eqref{eq:quadrel1}
by $(q^2y^4;q^2)_\infty=(qy^2,-qy^2;q)_\infty$ we see
that both sides become analytic in $y$.
By analytic continuation \eqref{eq:quadrel1} remains valid for
all $y^2\notin q^{-\N}$.

This proves the proposition for $z=\frac{q^{2-2k}}{t}$,
and by analytic continuation in $z$ the result follows.
\end{proof}

\begin{Rmk}\label{rmk:propquadrtransform}
(i) There are more choices possible for the
parameters in \eqref{eq:recrelfrombigqJacobi} to match
the recurrence \eqref{eq:recrelfrombigqJacobi} to
the recurrence
\eqref{eq:quadrrecrel} for the even and
odd degree $f_k$'s. All other possible choices lead to the
same result, In particular, going over the proof for the
odd degree $f_k$'s leads to \eqref{eq:quadrel1} with $t$
replaced by $qt$.

(ii) In \cite{KoelSCA} the spectral analysis of
the operator arising from the recurrence relation
\eqref{eq:recrelfrombigqJacobi} has been studied
on a suitable Hilbert space larger than $\lt$.
For the values of $a$ and $t$ as considered in
\S\S \ref{sec:specdecL1},
\ref{sec:specdecL2} this operator, say $S$, is, up to
a shift by a constant,
the square of $(L,{\mathcal D})$, so that
Theorems \ref{thm:case1spectraldecomp},
\ref{thm:case2spectraldecomp} also gives the spectral
decomposition of $S$. This shows that for the choices
of the parameters in \eqref{eq:recrelfrombigqJacobi}
as in the proof of Proposition \ref{prop:quadrtransform}
the spectral decomposition is explicit, cf. the remarks
on p.~193 and p.~200 of \cite{KoelSCA}.
Note that the parameters of
\eqref{eq:recrelfrombigqJacobi} used here
are not contained in the parameter set considered
in \cite{KoelSCA}.

(iii) Using Proposition \ref{prop:quadrtransform}
and the connection coefficients of Lemma \ref{lem:relsol}
and \cite[Prop.~4.4, Prop.~5.5]{KoelSCA} we can rewrite
any solution of \eqref{eq:diffequnnorm} in terms of
solutions of \eqref{eq:recrelfrombigqJacobi} in base $q^2$
with $(a,b,c,x)$ replaced by $(a,-1,-q,-t/q)$ or $(a,-1,-q,-t)$.
Apart from the case discussed in Proposition
\ref{prop:quadrtransform} this only gives more-term
transformations. As an example we give the expression for
the Ismail-Zhang $q$-analogue of the exponential
$\E_q$ defined in \eqref{eq:defcurlyE}.
Using Lemma \ref{lem:relsol} and
\eqref{eq:thetaidentity}, or \cite[(4.3.2)]{GaspR},
and Proposition \ref{prop:quadrtransform} in base $q^\hf$
gives, $z=\hf(y+y^{-1})$,
\begin{equation*}
\begin{split}
&\frac{(q^{\frac14}/y,-q^{\frac14}/y,-q^{\frac14}ty,
-q^{\frac14}/ty;q^\hf)_\infty(-q/t,-q^{\frac54}y^2/t;q)_\infty}
{(-q^\hf,y^{-2},-q^\hf/t;q^\hf)_\infty(q^{\frac34}/t;q^\hf)_\infty
(qt^2;q^2)_\infty}
\,_3\vp_2\left(
\genfrac{.}{.}{0pt}{}{q^{\frac34}y^2,-y^2,-q^\hf y^2}
{-q^{\frac54}y^2/t,qy^4};q,-\frac{q}{t}\right)\\ &\qquad\qquad+
\text{\rm Idem}(y\leftrightarrow y^{-1}) = \E_q(z;t)
\end{split}
\end{equation*}
\end{Rmk}


\bibliographystyle{amsplain}

\end{document}